\documentclass{article}
\usepackage{tikz}
\usepackage{hyperref}
\hypersetup{backref=true,       
	pagebackref=true,               
	hyperindex=true,                
	colorlinks=true,                
	breaklinks=true,                
	urlcolor= black,                
	linkcolor= green,                
	bookmarks=true,                 
	bookmarksopen=false,
	citecolor=black,
	linkcolor=blue,
	filecolor=black,
	citecolor=green,
	linkbordercolor=blue
}
\usepackage[all]{xy}
\usepackage{array,booktabs,arydshln,xcolor}
\usepackage[utf8]{inputenc}
\usepackage{lipsum}
\usepackage{ucs}
\usepackage{graphicx}
\usepackage[toc]{appendix}
\usepackage{pgf,pgfarrows,pgfnodes,pgfautomata,pgfheaps,pgfshade}
\usepackage{rotating}
\usepackage{multirow}
\usepackage{epstopdf}
\usepackage{mathtools}
\usepackage{amsmath,amssymb}
\usepackage{enumitem}
\usepackage{amsthm}
\usepackage{amsfonts}
\usepackage{authblk}
\usepackage{listings}
\usepackage{slashbox}
\usepackage{xcolor,framed,lipsum}

\usepackage{tikz,xcolor}
\usepackage{tikz-cd}
\usetikzlibrary{shapes.geometric,arrows,positioning,fit,calc,}
\usetikzlibrary{arrows}
\usepackage{imakeidx}
\usepackage{indentfirst}
\usepackage[left=3.5cm, right=3.5cm, top=2.5cm]{geometry}
\usepackage{array}
\usepackage{caption}

\newtheorem{theorem}{Theorem}[section]
\newtheorem*{theorem*}{Theorem}
\newtheorem{lemma}[theorem]{Lemma}
\newtheorem{cor}[theorem]{Corollary}
\newtheorem{prop}[theorem]{Proposition}

\theoremstyle{definition}
\newtheorem{definition}[theorem]{Definition}
\newtheorem{conc}[theorem]{Conclusion}

\newtheorem{notation}[theorem]{Notation}
\newtheorem{remark}[theorem]{Remark}

\newcommand\BibTeX{{\rmfamily B\kern-.05em \textsc{i\kern-.025em b}\kern-.08em
		T\kern-.1667em\lower.7ex\hbox{E}\kern-.125emX}}

\newcommand{\bigslant}[2]{{\raisebox{.2em}{$#1$}\left/\raisebox{-.2em}{$#2$}\right.}}

\title{\textbf{\textsc{SIMPLE SEMIGROUPS IN FINITE CATEGORIES\thanks{This project received funding from the European Research Council (ERC) under the European Union's Horizon 2020 research and innovation program (grant agreement No. 670624).\\
This work has been also supported by the French government, through the 3IA C\^ote d’Azur Investments in the Future project managed by the National Research Agency (ANR) with the reference number ANR-19-P3IA-0002.}}}}
\author[1,2]{Najwa Ghannoum \thanks{najwa.ghannoum@univ-cotedazur.fr}}
\author[1]{Carlos Simpson \thanks{carlos.simpson@univ-cotedazur.fr}}
\affil[1]{\small{Laboratoire J. A. Dieudonn\'e, CNRS, and Universit\'e C\^ote d'Azur}}
\affil[2]{Laboratoire LAMA, Universit\'e Libanaise}

\date{}
\begin{document}
	\maketitle
	\begin{abstract}
		In this paper, we classify finite categories with two objects such that one of the endomorphism monoids is a group. We prove that having a group on one side affects the structure of the other endomorphism monoid, and we prove that it is going to contain a simple semigroup. We also prove the other direction, that if we have a Rees matrix semigroup we can construct a category with two objects such that one of the objects is a group.\\
		
	\end{abstract}
	
	\section{Introduction}
	Studying finite categories mainly depends on studying the algebraic nature of their endomorphism monoids. In \cite{grouplike}, we classified finite categories with two objects such that their endomorphism monoids are grouplike. We studied the structure of such categories when the monoids contain a group. In this paper, we discuss categories such that one of the endomorphism monoids is a group, and in this case the group structure on one side affects the structure of the second monoid on the other side. In particular, the structure of the second monoid is going to contain a simple semigroup which is isomorphic to a Rees matrix semigroup with the same group \cite{JEP}.
	
	The notion of categories associated to matrices was initially introduced by T. Leinster and C. Berger in \cite{berger}, and was further discussed in more details by S. Allouch and C. Simpson in \cite{samer1}. Categories are associated to matrices with positive integers, such that the size of the matrix corresponds to the number of objects and the coefficients correspond to the number of morphisms between each two objects. In our work, we represent categories as matrices in terms of their monoids and bimodules \cite{grouplike}.
	
	\begin{definition}
		A bimodule is a set with actions on the left and the right of the respective monoids, such that the actions commute i.e. $(g \cdot x) \cdot h = g \cdot (x \cdot h)$. It can be seen as a category such that one of the sets of morphisms is empty, it's called an \textit{upper triangulated category}.
	\end{definition}
	
	\begin{definition}
		Let $\textsf{C}$ be a finite category. Then we get a matrix where the diagonal entries are monoids and the off diagonals entries are bimodules. This matrix is called \textit{algebraic matrix}.
	\end{definition}

\begin{definition}
	A category is called \textit{reduced} if it does not have any isomorphic distinct objects.
\end{definition}
	
	\begin{definition}
		A semigroup $S$ is \textit{simple} if its only ideals are $\emptyset $ and $S$. The Rees-Sushkevich structure theorem says that a simple semigroup is isomorphic to a \textit{Rees matrix semigroup}. In the following, we proceed without using this result, but some of the lemmas and constructions can also be viewed using this structure.
	\end{definition}
	
	\begin{remark}
		In our work we exclude the case where the ideals are empty. This means that whenever we say an ideal, then it's not empty, unless needed.
	\end{remark}

There are two directions:
	\begin{theorem} [\textbf{Monoids to categories}, Theorem \ref{montocatstheorem}] \label{theorem1} Let $A$ be a finite monoid, and not a group, and let $S$ be its minimal ideal. Then $S$ is a simple semigroup, and we can construct a category $\textsf{C}'$ with two objects associated to
		\begin{equation}\label{mat1} \tag{M(1)}
			\left(\begin{array}{cc} 
				S^{*1} & L \\
				R & G
			\end{array}\right) ~; ~ S^{*1} = I \cup \{1\} 
		\end{equation}
		where $L$ and $R$ are left and right minimal ideals of $S$ respectively and $G = L \cap R$ is a group. And
		$$
		|S| = \frac{|L|\cdot |R|}{|G|}.
		$$
		Then there is an inclusion $\textsf{C}' \subseteq \textsf{C}$ where $\textsf{C}$ is a category associated to the matrix
		\begin{equation}\label{mat2} \tag{M(2)}
			\left(\begin{array}{cc} 
				A & L \\
				R & G
			\end{array}\right).
		\end{equation}
	\end{theorem}
	
	\begin{remark} \label{remarkgroup}
	If $A$ is a group, then $S = A$ and in this case the identity of $S$ is the identity of $A$, and we obtain a groupoid \cite{grouplike}. In particular, we get $A \simeq G$ and $|L| = |R| = |G|$. Therefore, for this reason, we generally suppose in the following that $A$ is not a group.
	\end{remark}
%

	\begin{theorem} [\textbf{Categories to monoids}]  \label{goaltheorem} Let $\textsf{C}$ be a category associated to
		\begin{equation*}
			\left(\begin{array}{cc} \label{mat3} \tag{M(3)}
				A & U \\
				V & G
			\end{array}\right)
		\end{equation*}
		where $G$ is a group and $U,V$ not empty. Then $UV = S$ is the minimal ideal of $A$ and the construction of the category associated to (\ref{mat3}) is isomorphic to the construction of the category associated to (\ref{mat2}). And if $A$ is not a group then we get
		$$
			|A| \geq \frac{|L|\cdot |R|}{|G|} + 1.
		$$
	\end{theorem}
	
		\section{Categories with groups endomorphism monoids}
		In this section, we talk about categories with two objects where one of the endomorphism monoids is a group. The existence of a group on one of the diagonals of the matrix of a category imposes some properties over the sets of morphisms. These properties eventually characterize the nature of the other monoid.
		
		\begin{lemma} \label{freeaction}
			Let $\textsf{C}$ be a category associated to the algebraic matrix
			$$
			\left(\begin{array}{cc}
			A & L \\
			R & G
			\end{array}\right)
			$$
			such that $L$ and $R$ are not empty and $G$ is a group. Then $G$ acts freely on $L$ and $R$. Therefore $|L|$ and $|R|$ are multiples of $|G|$.
			\begin{proof}
			Let $g_1, g_2 \in G$ and $x \in L$ such that
			$$
			x \cdot g_1 = x \cdot g_2
			$$
			then 
			$$
			y \cdot x \cdot g_1 = y \cdot x \cdot g_2 ~;~y\in R.
			$$
		
		 We have that 
		 $$
		 y \cdot x = z \in G
		 $$
		 then
		  $$z^{-1}\cdot z \cdot g_1 = z^{-1} \cdot z \cdot g_2
		  $$
		 then $g_1=g_2$.
		
			Same for any $y\in R$.	
			\end{proof}
		\end{lemma}
	
	\begin{cor} \label{freeactioncor}
		Let $\textsf{C}$ be a category associated to the algebraic matrix
		$$
		\left(\begin{array}{cc}
			A & L \\
			R & G
		\end{array}\right)
		$$
		such that $L$ and $R$ are not empty and $G$ is a group. Then $G$ acts freely on $L \times R$.
		\begin{proof}
			Let $g_1,g_2 \in G$ and let $(x,y) \in L \times R$ such that
			\begin{eqnarray*}
			g_1 \cdot (x,y) &=& g_2 \cdot (x,y) \\
			(x \cdot g_1, g_1 \cdot y) &=& (x \cdot g_2, g_2 \cdot y) 
			\end{eqnarray*}
		then
		$$
		x \cdot g_1 = x \cdot g_2 ~\textrm{and}~ g_1 \cdot y = g_2 \cdot y
		$$
		and by the free action of $G$ on $L$ and $R$ we get
		$$
		g_1 = g_2.
		$$
		\end{proof}
	\end{cor}

\begin{remark}
	The condition that $L$ and $R$ are not empty is important to obtain the results. Therefore, in the following we always assume that the bimodules $L$ and $R$ are not empty.
\end{remark}
		
		Let $B$ be a monoid, and let $L$ be a right $B$-module and let $R$ be a left $B$-module. Then there exists a set
		$$
		L \otimes_B R := L \times R / \sim_B
		$$
		such that
		$$
		(xb,y) \sim_B (x,by) ~ \forall x \in L, b\in B, y\in R.
		$$
		$\sim_B$ is an equivalence relation.
		
		
		Then there's a map
		$$
		L \times R \rightarrow L \otimes_B R := L \times R /\sim_B \rightarrow LR~;~ (x,y) \mapsto x\otimes y \mapsto x\cdot y.
		$$
		
		\begin{lemma}\label{assoc}
			Let $A,B,A'$ be monoids.
			\begin{enumerate}
				\item If $L$ is an $(A,B)$-bimodule then $L \otimes_B R$ is a left $A$-module.
				\item If $R$ is a $(B,A')$-bimodule then $L \otimes_B R$ is a right $A'$-module.
				\item If $L,R$ are $(A,B)$ and $(B,A)$ bimodules, then $L \otimes_B R$ is an $(A,A)$-bimodule.
			\end{enumerate}
		\end{lemma}
		
		\begin{lemma}
			If $B=G$ is a group, then $G$ acts on $L \times R$ and
			$$
			(x,y) \sim (x',y') \iff \exists~ g\in G ~\textrm{such that}~ x' = xg^{-1} ~\textrm{and}~ y' = gy.
			$$
			In this case
			$$
			g \cdot (x,y) := (xg^{-1}, gy)~\textrm{and}~ L \otimes_G R = L \times R /G.
			$$
		\end{lemma}
		
		\begin{remark}
			If $B=G$ is a group, we can transform a left action to a right action and vice-versa by multiplying with inverses.
		\end{remark}
		
		\begin{prop}\label{determinantprop}
			Let $\textsf{C}$ be a finite category associated to the matrix
			\begin{center}
				$\left( \begin{array}{cc}
					A & L \\
					R & G
				\end{array} \right)$.
			\end{center}
			Consider the quotient $L \times R /G$ with the equivalence relation
			$$
			(x,y) \sim (x',y') \iff \exists~ g\in G ~\textrm{such that}~ x' = xg^{-1} ~\textrm{and}~ y' = gy.
			$$
			Then the function 
			$$
			f : L \times R / G \rightarrow LR~;~ f(x,y) = xy
			$$
			is bijective. In this case $|LR| = |L \times R / G| = \frac{|L||R|}{|G|}$. 
			
			\begin{proof}
				\begin{enumerate}
					\item \textbf{$f$ is injective:} Let $(x,y), (x',y') \in L\times R /G$ such that $f(x,y) = f(x',y')$ then $xy = x'y'$. There exist $c \in L$ and $d\in R$ such that $yc = 1_G$ and $dx = 1_G$. This gives us
					$$
					x = x \cdot (yc) = x'(y'c)
					$$
					and
					$$y = (dx) \cdot y = (dx')y'
					$$ where $y'c$ and $dx'$ are in $G$ and they are inverses, indeed
					$$
					(dx')(y'c) = d(x'y')c = d(xy)c = (dx)(yc) = 1_G.
					$$
					Then $(x,y) \sim (x',y')$.
					\item \textbf{$f$ is surjective:} Evident.
					\item Since $G$ acts freely on $L$ and on $R$ (Lemma \ref{freeaction}) then it also acts freely on $L \times R$ (Corollary \ref{freeactioncor}). Then
					$|L \times R / G| = \frac{|L||R|}{|G|}$.
				\end{enumerate}
				
			\end{proof}
		\end{prop}
	\begin{conc}
	Let $\textsf{C}$ be a category associated to the matrix
	$$
	\left(\begin{array}{cc}
	A & L \\
	R & G	
	\end{array}\right)
	$$
	where $G$ is a group. Then
	\begin{enumerate}
		\item $G$ acts freely on $L$ and $R$.
		\item $|L|$ and $|R|$ are multiples of $|G|$.
		\item If $A$ is not a group then $|A| \geq \frac{|L||R|}{|G|}+1$.
	\end{enumerate}
	\end{conc}
		
		\section{From simple semigroups to categories} \label{simpletocats}
		We aim in this section to construct a finite category with two objects where one of the endomorphism monoids is a group. We start with a simple semigroup $S$, to which we add an identity element, then we take the minimal left and right ideals of $S$, denote them by $L$ and $R$ respectively. The intersection of $L$ and $R$ is a group, which will be the second endomorphism monoid of the category.
		
		\begin{definition}
			Given a monoid (or semigroup) $S$, a \textit{left ideal} in $S$ is a subset $A$ of $S$ such that $SA$ is contained in $A$. Similarly, a \textit{right ideal} is a subset $A$ such that $AS$ is contained in $A$. Finally, a \textit{two-sided ideal}, or simply \textit{ideal}, in $S$ is a subset $A$ that is both a left and a right ideal.
		\end{definition}
		
		\begin{prop}
			Let $S$ be a simple semigroup, $L$ be a minimal left ideal and $R$ be minimal right ideal of $S$. Then $LR = S$.
			\begin{proof}
				We have $L \subset S$ and $R \subset S$ then $LR \subset S$. In addition, $LR$ is a two sided ideal of $S$. Indeed, let $x \in S$ and $ab \in LR$ such that $a \in L$ and $b \in R$, then 
				$$
				x (ab) = \underbrace{(xa)}_{\in L}b \in LR
				$$
				and
				$$
				(ab)x = a\underbrace{(bx)}_{\in R} \in LR.
				$$
				Hence $LR \subset S$ is a two sided ideal, but $S$ is simple, then $LR = S$.
			\end{proof}
		\end{prop}
		
		\begin{theorem}\label{group}
			Let $S$ be a simple semigroup, $L$ be a minimal left ideal and $R$ be minimal right ideal of $S$. Then:
			\begin{enumerate}
				\item $G= L\cap R$ is a group.
				\item $RL = L\cap R$.
			\end{enumerate}
			\begin{proof}\
				For part (1), it is sufficient to prove that G is right and left cancellative.\\
				Let $z\in L \cap R$ and $zR$ be the right ideal contained in $R$, $zR = \{zr_1,...,zr_n\}$, let $r_i,r_j$ such that $zr_i=zr_j$, if $r_i\neq r_j$ then $zR \subsetneq R$. This is a contradiction with the minimality of $R$.\\
				Let $Lz$ be the left ideal contained in $L$, $Lz = \{l_1z,...,l_mz\}$, let $l_i,l_j$ such that $l_iz=l_jz$, if $l_i\neq l_j$ then $Lz \subsetneq L$. This is a contradiction with the minimality of $L$.\\
				Hence, $G$ is right and left cancellative. Proving now the existing of an identity element.\\
				Since G is right cancellative, then the map $x\mapsto x\cdot z$ is injective, as $G$ is finite then it is also surjective. Then for all $w$, there exists $e; ez=w$. Take $w=z$.\\
				Also, $G$ is left cancellative, then the map $x\mapsto z\cdot x$ is bijective, and for all $y$ there exists $e'$ such that $ze'=y$. Take $y=z$.\\
				Then:\\
				- $ey=e(zx)=(ez)x=zx=y$.\\
				- $we'=(ez)e'=ez=w$.\\
				Hence the existence of an identity element.\\
				
				For part (2), we always have $RL \subseteq L\cap R$.\\
				Let $x\in L\cap R=G$, $x=e \cdot x \in RL$ where $e$ is the identity of $G$.
			\end{proof}
		\end{theorem}


\begin{theorem} \label{simpletocatstheorem}
	Let $S$ be a simple semigroup. Let $L$ and $R$ be minimal left and right ideals of $S$ respectively. Then we can construct a category $\textsf{C}$ associated to the matrix
	$$
	\left( \begin{array}{cc}
S^{*1} & L \\
R & G
	\end{array} \right)~;~ S^{*1} = S \cup \{1\}
	$$
	such that $G = RL$ is a group and $LR = S$ where
	\begin{itemize}
		\item $|L|$ and $|R|$ are multiples of $|G|$.
		\item $|S| = |LR| = \frac{|L||R|}{|G|}$. 
	\end{itemize}
	\begin{proof}
		We want to construct a category with two objects starting with a semigroup. The construction we use is called \textit{Karoubi envelope} or \textit{idempotent completion}. It was first introduced in 1963 by M. Artin, A. Grothendieck and J.L. Verdier \cite{artin269theorie}, then the definition also appears in \cite{borceux1986cauchy} and \cite{amarilli2021locality}. The idea is that if $\textsf{C}$ is a category, then its idempotent envelope $\textsf{C}^{idem}$ (also denoted $\tilde{\textsf{C}}$) is the category whose objects are pairs $(X,p)$ where $X \in Ob(\textsf{C})$ and $p : X \rightarrow X$ is an idempotent. The set of morphisms from $X$ to $X'$ is
		$$
		\textsf{C}^{idem}((X,p),(X',p')) = p \cdot \textsf{C}(X,X')\cdot p \subseteq \textsf{C}(X,X').
		$$
		
		Then in the case of a monoid $A$, we view it as a category with one object $\{\ast\}$ and $A$ is the set of endomorphisms of $\{\ast\}$, we denote the category by $(\{\ast\}, A)$. The Karoubi envelope will be the category
		$$
		Idem(A) := (\{\ast\}, A)^{idem}
		$$
		where the objects of $Idem(A)$ are the idempotents of $A$.
		
		If we apply this construction to our work here, we get the following: if $S$ is a simple semigroup, $L$ minimal left ideal of $S$, $R$ minimal right ideal of $S$ and $G$ the group obtained by $L \cap R$. Then we take the two idempotents $e_1 = 1_{S^{*1}}$ and $e_2 = 1_G$, and we want to construct a category $Idem_{e_1,e_2}(S^{*1})$, the sets of morphisms are the following
		\begin{center}
			$S^{*1} = e_1 \cdot S^{*1} \cdot e_1 ~~~~~~~~~~~~ L = e_1 \cdot S^{*1} \cdot e_2$ \\
			$R = e_2 \cdot S^{*1} \cdot e_1 ~~~~~~~~~~~~ G = e_2 \cdot S^{*1} \cdot e_2$.
		\end{center}
	Note that $L$ and $R$ here are the same of Theorem \ref{group}, and they are left and right minimal ideals of $S$ respectively.
	
		Then $Idem_{e_1,e_2}(S^{*1})$ is a full subcategory of $Idem(S^{*1}) = \widetilde{S^{*1}}$ the Karoubi envelope of $S^{*1}$. The multiplication of morphisms is defined by the multiplication of elements of $S$.
	\end{proof}
\end{theorem}

\begin{conc}
In this section, we started with a simple semigroup $S$ and we proved that we can construct a category with two objects such that one of the objects is $S$ and the other one is a group.
\end{conc}
	
		\section{From finite monoids to categories}
		In this section, we study how we can construct a category with two objects such that one of the objects is a group by starting with a finite monoid instead of a simple semigroup.
	
		\begin{prop}
			Let $A$ be a finite monoid, then $A$ has a unique non-empty minimal ideal denoted by $S_0$.
			\begin{proof}
				Let $S = \{S_i \subset A \mid S_i ~ \textrm{ideal} ~ of A\} = \{S_1, S_2, \hdots , S_k\}$ be the set of all ideals of $M$, set $S_0 = \underset{S_i \in S}{\bigcap}S_i = S_1 \cap S_2 \cap \hdots \cap S_k$. $S_0$ is not empty because it contains at least $S_1 \cdot S_2 \cdot \hdots \cdot S_k$ and it is minimal because it doesn't contain any other ideal. $S_0$ is unique. Indeed, suppose that $S'$ is another minimal ideal, then $S_0 \subseteq S'$ because $S_0$ is the intersection of all ideals, but $S'$ is minimal then $S_0 = S'$.
			\end{proof}
		\end{prop}
		
			\begin{remark}
			Let $A$ be a finite monoid and $S_0$ be the minimal ideal of $A$, then $S_0$ is a sub-semigroup of $A$.
		\end{remark}
	
		\begin{lemma} \label{leftofM}
			Let $A$ be a finite monoid. If $M \subseteq A$ is a left ideal, then $M$ is a sub-monoid of $A$ and if $L \subseteq M$ is a minimal left ideal of $M$ then $L$ is a minimal left ideal of $A$.
			\begin{proof}
				We have $M \cdot L \subseteq L$ and $M \cdot L \subseteq M$ is a left ideal of $M$ contained in $L$, by the minimal of $L$ in $M$ we have $M \cdot L = L$, i.e. for all $l \in L$ there exists $m \in M, l' \in L$ such that $l = ml'$.
				
				Now let $a \in A$ and $l \in L$ then
				\begin{eqnarray*}
					a \cdot l &=& a \cdot (ml')\\
							  &=& (am)\cdot l' \\
							  &\in& M \cdot L~\textrm{(because $M$ is left ideal of $A$)} \\
							  &\in& L~\textrm{(because $L$ is left ideal of $M$)}.
				\end{eqnarray*}
			Therefore, $L$ is a left ideal of $A$.
			
			$M$ is minimal in $A$. Indeed, suppose that $L'$ is a left ideal of $A$ such that $L' \subseteq L$, then
			$$
			M \cdot L' \subseteq L' \subseteq L
			$$
			is a left ideal of $M$, but $L$ is minimal in $M$, then $M\cdot L' = L$. Therefore, $L' = L$.
			\end{proof}
		\end{lemma}
	
	\begin{lemma} \label{rightofM}
		Let $A$ be a finite monoid. If $M \subseteq A$ is a right ideal, then $M$ is a sub-monoid of $A$ and if $R \subseteq M$ is a minimal right ideal of $M$ then $R$ is a minimal right ideal of $A$.
		\begin{proof}
			Similar to the proof of Lemma \ref{leftofM}.
		\end{proof}
	\end{lemma}
		
		\begin{lemma} \label{minofAminofS}
			Let $A$ be a finite monoid and let $S_0$ be its minimal idea. Then $L$ (resp. $R$) is a minimal left (resp. minimal right) ideal of $A$ if and only if $L$ (resp. $R$) is a minimal left (resp. minimal right) ideal of $S_0$. 
			\begin{proof}
				We have $S_0\cdot L \subset S_0$ and $S_0 \cdot L \subset L$ is a left ideal of $A$, by the minimality of $L$ in $A$ we obtain that $S_0\cdot L = L$. Then $L \subset S_0$. Also, $L$ is a minimal ideal of $S_0$. Indeed, if $L' \subset L$ is a left ideal of $S_0$ then $S_0 \cdot L' \subset L' \subset L$, but $S_0 \cdot L'$ is a left ideal of $A$ and $L$ is minimal of $A$ then $S_0 \cdot L' = L$ hence $L' = L$ and $L$ is a minimal left ideal of $S_0$.
				
				For the other direction, use Lemma \ref{leftofM} and Lemma \ref{rightofM}.
			\end{proof}
		\end{lemma}
		
		\begin{lemma}
			Let $A$ be a finite monoid and $S_0$ the minimal ideal of $A$, then $S_0$ is a simple semigroup. In particular, $S_0$ has the structure of a Rees matrix semigroup.
			\begin{proof}
				If $J \subset S_0$ is an ideal of $S_0$ then $S_0\cdot J \cdot S_0 \subset J \subset S_0$ and $S_0 \cdot J \cdot S_0$ is an ideal of $A$, then by the minimality of $S_0$ in $A$, we obtain $S_0 \cdot J \cdot S_0 = S_0$ hence $J = S_0$ and $S_0$ is simple.
			\end{proof}
		\end{lemma}
		
		\begin{theorem} \label{montocatstheorem}
			Let $A$ be a finite monoid that's not a group, and let $S_0$ be its minimal ideal. Then by Section \ref{simpletocats}, we obtain a category $\textsf{C}'$ associated to the matrix
			$$
			\left(\begin{array}{cc}
				S_0^{*1} & L \\
				R & G
			\end{array}\right)
			$$
			where $L$ and $R$ are minimal left and right ideals of $S_0$ respectively, and $G$ is a group.
			
			Then there exists a category $\textsf{C}$ associated to the matrix
			$$
			\left(\begin{array}{cc}
				A & L \\
				R & G
			\end{array}\right)
			$$
			such that $\textsf{C}' \subseteq \textsf{C}$.
			\begin{proof}
				We use the same construction of Theorem \ref{simpletocatstheorem}. Then if we apply this construction to our work here, we get the following: if $A$ is a monoid, $L$ minimal left ideal of $S$, $R$ minimal right ideal of $A$ and $G$ the group obtained by $L \cap R$. Then we take the two idempotents $e_1 = 1_{A}$ and $e_2 = 1_G$, and we want to construct a category $Idem_{e_1,e_2}(A)$, the sets of morphisms are the following
				\begin{center}
					$A = e_1 \cdot A \cdot e_1 ~~~~~~~~~~~~ L = e_1 \cdot A \cdot e_2$ \\
					$R = e_2 \cdot A \cdot e_1 ~~~~~~~~~~~~ G = e_2 \cdot A \cdot e_2$.
				\end{center}
				 $L$ is a left minimal ideal of $A$ and $R$ is a minimal right ideal of $A$.
				
				Then $Idem_{e_1,e_2}(A)$ is a full subcategory of $Idem(A) = \tilde{A}$ the Karoubi envelope of $A$. The multiplication of morphisms is defined by the multiplication of elements of $A$.
				
				Now since $A$ is not a group, then $1_A \notin S_0$ (Remark \ref{remarkgroup}), hence
				$$
				S_0^{*1} \rightarrow A
				$$
				is injective.
				
				 $L$ and $R$ are minimal left and right ideals of $S_0$ respectively. By Lemma \ref{minofAminofS}, $L$ and $R$ are also minimal left and right ideals of $A$ respectively, i.e. for all $x\in L, y\in R$ and $a\in A$, we have $a \cdot x \in L$ and $y \cdot a \in R$. By Section \ref{simpletocats} we have $LR = S_0$. Therefore $LR = S_0 \subset A$.
				
				In addition, for all $x \in L$ and $y \in R$, $x\cdot y \neq 1_A$, because $G$ is a group and in this case we get isomorphic objects, which gives us that $A$ is a group and this is a contradiction with the hypothesis.
			\end{proof}
		\end{theorem}
	\begin{conc}
		For any finite monoid $A$, we can construct a category with two objects such that one of the objects is a group. This result means that every finite monoid $A$ is connected to a group, which is the intersection of a minimal left and a minimal right ideal of $A$. Therefore, Theorem \ref{theorem1} is proved.
	\end{conc}
		
		\section{From categories to simple semigroups}
		In this section, we give a proof of the direction that starting with a category $\textsf{C}$ with two objects such that one of the objects is a group, then we can obtain a simple semigroup with cardinality as shown in the previous sections.
		
		\begin{theorem}\label{determinant}
			Let $\textsf{C}$ be a finite category associated to the matrix
			\begin{center}
				$\left( \begin{array}{cc}
					A & U \\
					V & G
				\end{array} \right)$
			\end{center}
			where $G$ is a group and $A$ is not a group, then $S = UV \subseteq A$ is a simple ideal such that $|S| = \frac{|U||V|}{|G|}$ and there exists a sub-category $\textsf{C}'$ of $\textsf{C}$ associated to the matrix
			\begin{equation}\label{rees2}
				\left( \begin{array}{cc}
					S^{*1} & U \\
					V & G
				\end{array} \right)
			\end{equation}
			where $S^{*1} = S \cup \{1\}$ and $S$ is a simple semigroup.
			\begin{proof}
				$\textsf{C}'$ is a sub-category of $\textsf{C}$:
				\begin{itemize}
					\item $Ob(\textsf{C}') = Ob(\textsf{C})$.
					\item Morphisms $= \{S^{*1}, U, V, G\}$.	
				\end{itemize}
				
				Suppose $S' \subset S $ is a two sided ideal, $S \subset S'$?
				
				Let $x \in U, y\in V$ and $a \in S'$. We have $y \cdot a \cdot x \in G $, then there exists $g\in G$ such that
				\begin{eqnarray*}
					g \cdot (y \cdot a \cdot x) &=& 1_G
				\end{eqnarray*}
					then
					\begin{eqnarray*}
					x \cdot y &=& x \cdot 1_G \cdot y \\
					&=& x \cdot (g\cdot y\cdot a \cdot x) \cdot y \\
					&=& \underbrace{(x \cdot g \cdot y)}_{\in S} \cdot \underbrace{a}_{\in S'} \cdot \underbrace{(x \cdot y)}_{\in S} \in S'.
				\end{eqnarray*}
			\end{proof}
		\end{theorem}
	
	\begin{definition}
		Let $\textsf{C}$ be a finite category associated to the matrix
		\begin{center}
			$\left( \begin{array}{cc}
				A & U \\
				V & G
			\end{array} \right)$
		\end{center}
		where $G$ is a group. Let $x \in U$ and $y\in V$, define $L_y$ and $R_x$ in the following way
		$$
		L_y := U\cdot y \subseteq A
		$$
		and
		$$
		R_x := x \cdot V \subseteq A.
		$$
	\end{definition}

\begin{lemma}
$L_y$ is a minimal left ideal of $A$ and $R_x$ is a minimal right ideal of $A$.
\begin{proof}
It is clear that $L_y$ is a left ideal of $A$. Proving it's minimal. Suppose $L' \subseteq L_y$ a left ideal of $A$. Let $u \in U, a\in L'$, we have $y\cdot a \cdot u \in G$ then there exists $g\in G$ such that
$$
g \cdot (y \cdot a \cdot u) = 1_G
$$
then
	\begin{eqnarray*}
	u \cdot y &=& u \cdot 1_G \cdot y \\
	&=& u \cdot (g\cdot y\cdot a \cdot u) \cdot y \\
	&=& \underbrace{(u \cdot g \cdot y)}_{\in L_y} \cdot \underbrace{a}_{\in L'} \cdot \underbrace{(u \cdot y)}_{\in L_y} \in L'.
\end{eqnarray*}
Hence $L_y$ is minimal. Similarly for $R_x$.
\end{proof}
\end{lemma}

\begin{remark}
	For all $g\in G$, if we replace $y$ by $gy$ we obtain
	$$
	L_{gy} = U \cdot g \cdot y = U \cdot y = L_y~(\textrm{because}~U \cdot g = U).
	$$
\end{remark}

	\begin{prop} \label{bijection}
	Let $\textsf{C}$ be a category associated to the matrix
	$$
	\left(\begin{array}{cc}
		A & U \\
		V & G
	\end{array}\right)
	$$
	Let $\mathcal{L}(A) = \{L \subseteq A \mid L~ \textrm{is a minimal left ideal of $A$}\}$ and $\mathcal{R}(A) = \{R \subseteq A \mid R ~\textrm{is a minimal right ideal of $A$}\}$.
	
	Define the sets 
	$$
	\bigslant{\{y\}}{y \sim gy} :=~ _{G}\backslash^{V} ~~~\textrm{and}~~~ \bigslant{\{x\}}{x \sim xg} := \bigslant{_U}{_G}.
	$$
	Then
	$$
	_{G}\backslash^{V} = \mathcal{L}(A) ~~\textrm{and}~~\bigslant{_U}{_G} = \mathcal{R}(A).
	$$
	This means that for every minimal left ideal $L$ of $A$, there exists $y \in V$ such that $L = L_y$. Similarly, for every minimal right ideal $R$ of $A$, there exists $x \in U$ such that $R = R_x$.
\end{prop}

\begin{lemma}
	Let $\textsf{C}$ be a finite category associated to the matrix
	\begin{center}
		$\left( \begin{array}{cc}
			A & U \\
			V & G
		\end{array} \right)$
	\end{center}
	where $G$ is a group and $S= UV$ a simple ideal of $A$. For $x \in U$ and $y \in V$ such that $y\cdot x = 1_G$, we have
	$$
	L_y R_x = U \cdot y \cdot x \cdot V = UV = S.
	$$
\end{lemma}	

	\begin{prop}
	Let $A$ be a finite monoid that is not a group, let $L$ be a minimal left ideal and $R$ be a minimal right ideal of $A$, let $x,y \in A$, then
	\begin{center}
		$(xR)\cap (Ly) = xRLy$.
	\end{center}
	\begin{proof}
		Theorem \ref{group}.
	\end{proof}
\end{prop}

\begin{notation}
	$(xR)\cap (Ly) = xRLy= G_{xy}$.	
\end{notation}

\begin{theorem} \label{catstosimpletheorem}
		Let $\textsf{C}$ be a finite category associated to the matrix
	\begin{center}
		$\left( \begin{array}{cc}
			A & U \\
			V & G
		\end{array} \right)$
	\end{center}
	where $G$ is a group. Let $x \in U, y\in V$ such that $y\cdot x = 1_G$, then by the construction described in Theorem \ref{montocatstheorem} there exists a category $\textsf{C}'$ associated to the matrix
	\begin{center}
		$\left( \begin{array}{cc}
			A & L_y \\
			R_x & G_{xy}
		\end{array} \right)$
	\end{center}
such that $L_y R_x = S$ the simple ideal of $A$ where $L_y$ is a minimal left ideal of $A$ and $R_x$ is a minimal right ideal of $A$ and $G_{xy} = x \cdot G \cdot y$ is a group. Then $\textsf{C}' \simeq \textsf{C}$.
	\begin{proof}
	Consider the maps
	$$
	\phi : G \rightarrow G_{xy}~;~ g \mapsto xgy
	$$
	
	$$
	\psi : U \rightarrow L_y ~;~ u \mapsto uy
	$$
	and
	$$
	\psi' : V \rightarrow R_x ~;~ v \mapsto xv.
	$$
	$\phi, \psi$ and $\psi'$ provide the isomorphisms needed to prove the theorem.
\end{proof}
\end{theorem}

\begin{prop}
	The category obtained is unique.
	\begin{proof}
		Let $A$ be a monoid, choose $L,R$ such that $G = RL$. Choose another $L',R'$ such that $G' = R'L'$.
		
		Since we have $G= RL$ then we can choose $L_y,R_x \subseteq A$ such that
		$$
		\textsf{C}(A,L,R,G) \simeq \textsf{C}(A,L_y,R_x,G_{xy}).
		$$
		Take $L_y = L'$ and $R_x = R'$, then
		$$
		\textsf{C}(A,L',R',G) = \textsf{C}(A,L_y,R_x,G_{xy}).
		$$
		Hence $\textsf{C}(A,L,R,G) \simeq \textsf{C}(A,L',R',G')$.
	\end{proof}
\end{prop}

\begin{conc}
	Finally in this section, we prove that if we have a category $\textsf{C}$ associated to the matrix
	\begin{center}
		$\left( \begin{array}{cc}
			A & U \\
			V & G
		\end{array} \right)$
	\end{center}
then under a choice of an element $x\in U$ and an element $y\in V$, we obtain a left minimal ideal $L_y$ and a right minimal ideal $R_x$ of $A$ such that $L_yR_x$ is a simple semigroup and $R_xL_y$ is a group. This data gives us a unique category $\textsf{C}'$ associated to the matrix
\begin{center}
	$\left( \begin{array}{cc}
		A & L_y \\
		R_x & G_{xy}
	\end{array} \right)$
\end{center}
such that $\textsf{C}' \simeq \textsf{C}$ and $|L_yR_x| = \frac{|L_y||R_x|}{|G_{xy}|}$. If $A$ is not a group, then we get
$$
|A| = \frac{|U|\cdot|V|}{|G|} +1.
$$
Therefore, Theorem \ref{goaltheorem} is proved.
\end{conc}

	\section{Connectivity of monoids}
	\begin{definition}
		We say that two monoids $A$ and $B$ are \textit{connected} if there exists a category with two objects such that its endomorphism monoids are $A$ and $B$.
	\end{definition}
	
	\begin{definition}
		We say that a monoid $A$ \textit{has a group} $G$ if its minimal ideal is connected to $G$.
	\end{definition}

\begin{remark}
	Let $A$ be a finite monoid, by the discussion in the previous chapters, there is always a group connected to $A$ through the minimal ideal $S$ of $A$. The group is of the form $x \cdot G \cdot y$ where $x,y \in S$ and $G = RL$ where $L$ and $R$ are the minimal left and right ideals of $S$ (eventually $A$) respectively. In addition, $S$ is a simple semigroup, then by Rees-Sushkevich theorem, it is a Rees matrix semigroup, and the group $G$ is the group that is involved in the definition of the Rees matrix semigroup.
\end{remark}

	We would like to thank J\'er\'emie Marqu\`es for suggesting this result in the following lemma.
	\begin{lemma} [J\'er\'emie Marqu\`es] \label{transitive}
		If $A,B$ are connected, and $B,C$ are connected then $A,C$ are connected. 
		\begin{proof}
			$A,B$ are connected, then there exists a category $\textsf{C}_1$ associated to the algebraic matrix
			$$
			M_1 = \left(\begin{array}{cc}
				A & L \\
				R & B
			\end{array}\right)
			$$
			$B,C$ are connected, then there exists a category $\textsf{C}_2$ associated to the algebraic matrix
			$$
			M_2 = \left(\begin{array}{cc}
				B & L' \\
				R' & C
			\end{array}\right).
			$$
			Let
			$$
			M = \left(\begin{array}{cc}
				A & L \otimes_B L' \\
				R' \otimes_B R & C
			\end{array}\right)
			$$
			then $M$ is a matrix of a category $\textsf{C}$. Indeed,	$A$ acts to the left on $L$ and $C$ acts to the right on $L'$, then $A$ and $C$ act on $L \otimes_B L'$. Same for $R$ and $R'$. It remains to prove $(L \otimes_B L') \otimes_C (R' \otimes_B R) \rightarrow A$ and $(R' \otimes_B R) \otimes_A (L \otimes_B L') \rightarrow B$.
			\begin{eqnarray*}
				(L \otimes_B L') \otimes_C (R' \otimes_B R) &\simeq& L \otimes_B (L' \otimes_C R') \otimes_B R \\
				&\rightarrow& L \otimes_B B \otimes_B R\\
				&\simeq& L \otimes_B R \\
				&\rightarrow& A.
			\end{eqnarray*}
		
			We have that 
			\begin{equation} \label{tens1}
				L \otimes_B R \rightarrow A
			\end{equation}
			is an $(A,A)$-bimodule morphism,	and 
			\begin{equation} \label{tens2}
				R \otimes_A L \rightarrow B	
			\end{equation}
			is a $(B,B)$-bimodule morphism.
			
			To verify associativity, the morphisms (\ref{tens1}) and (\ref{tens2}) should satisfy the conditions that the following diagrams are commutative.
			\begin{center}
				\begin{tikzcd}
					L\otimes_B R \otimes_A L \arrow[r] \arrow[d] & A \otimes_A L \arrow[d] \\
					L \otimes_B B \arrow[r]                      & L                      
				\end{tikzcd} ~\textrm{and}~ \begin{tikzcd}
					R\otimes_A L \otimes_B R \arrow[r] \arrow[d] & B \otimes_B R \arrow[d] \\
					R \otimes_A A \arrow[r]                      & R                      
				\end{tikzcd}
			\end{center}
			Same for $(R' \otimes_B R) \otimes_A (L \otimes_B L')$.
		\end{proof}
	\end{lemma}

	\begin{theorem}
		Two monoids are connected if and only if they have the same group.
		\begin{proof}
			For the first direction, let $A$ be a monoid that has a group $G$, then the minimal ideal $S_0$ of $A$ is Rees matrix semigroup and connected to $G$.
			
			Similarly, Let $B$ be a monoid that has a group $H$, then the minimal ideal $J_0$ of $A$ is Rees matrix semigroup and connected to $H$.
			
			If $A$ is connected to $B$, then by transitivity (Lemma \ref{transitive}), $G$ and $H$ are connected, then $G$ and $H$ are isomorphic (Remark \ref{remarkgroup}).
			
			For the second direction, let $A$ and $B$ be two monoids that have the same group $G$, then $A$ is connected to $G$ and $B$ is connected to $G$, by transitivity (Lemma \ref{transitive}), $A$ and $B$ are connected.
		\end{proof}
	\end{theorem}

	\bibliographystyle{siam}
	\bibliography{biblio}

\end{document}